\documentclass[a4paper,twoside,reqno]{amsart}

\usepackage[pagewise]{lineno}%\linenumbers
\usepackage[utf8]{inputenc}

\usepackage{authblk}
\usepackage{hyperref}
\hypersetup{urlcolor=blue, citecolor=red}
\usepackage{amsmath,amsthm,amssymb,xcolor,bbm,mathrsfs}
\usepackage[english]{babel}
\usepackage{fancyhdr}
\newcommand{\norm}[1]{\|#1\|}
\newcommand{\n}{\hspace*{-6pt}}
\DeclareMathOperator{\dist}{dist}
\DeclareMathOperator{\diam}{diam}
\DeclareMathOperator{\diag}{diag}
\DeclareMathOperator{\unif}{uniform}

\makeatletter
\@namedef{subjclassname@1991}{2020 Mathematics Subject Classification}
\makeatother

\title{Mixed Hegselmann-Krause Dynamics\\ %\uppercase\expandafter{\romannumeral2}
\textemdash\footnotesize nondeterministic case}
\author{Hsin-Lun Li}

\date{}
\email{hsinlunl@asu.edu}

\theoremstyle{definition}
\newtheorem{theorem}{Theorem}
\newtheorem{definition}{Definition}
\newtheorem{lemma}{Lemma}
\newtheorem{corollary}{Corollary}
\newtheorem{property}{Property}
\newtheorem{example}{Example}

\begin{document}

\allowdisplaybreaks

\thispagestyle{firstpage}
\maketitle
\begin{center}
    Hsin-Lun Li\\
    School of Mathematical and Statistical Sciences,\\ Arizona State University, Tempe, AZ 85287, USA
\end{center}

\begin{abstract}
The original Hegselmann-Krause (HK) model is composed of a finite number of agents characterized by their opinion, a number in $[0,1]$. An agent updates its opinion via taking the average opinion of its neighbors whose opinion differs by at most $\epsilon$ for $\epsilon>0$ a confidence bound. An agent is absolutely stubborn if it does not change its opinion while update, and absolutely open-minded if its update is the average opinion of its neighbors. There are two types of HK models--the synchronous HK model and the asynchronous HK model. The paper is about a variant of the HK dynamics, called the mixed model, where each agent can choose its degree of stubbornness and mix its opinion with the average opinion of its neighbors at all times. The mixed model reduces to the synchronous HK model if all agents are absolutely open-minded all the time, and the asynchronous HK model if only one uniformly randomly selected agent is absolutely open-minded and the others are absolutely stubborn at all times. In \cite{mhk}, we discuss the mixed model deterministically. Point out some properties of the synchronous HK model, such as finite-time convergence, do not hold for the mixed model. In this topic, we study the mixed model nondeterministically. List some properties of the asynchronous model which do not hold for the mixed model. Then, study circumstances under which the asymptotic stability holds. 
\end{abstract}

\section{introduction}
The original Hegselmann-Krause (HK) model involves a finite number of agents characterized by their opinion, a number in $[0,1].$ Agent $i$ updates its opinion $x_i$ via taking the average opinion of its neighbors whose opinion differs by at most $\epsilon$ from $x_i$ for $\epsilon>0$ a confidence bound. There are two types of HK models--the synchronous HK model and the asynchronous HK model. For the synchronous HK model, all agents update their opinion at each time step, whereas for the asynchronous HK model, only one agent chosen uniformly at random updates its opinion at each time step. The mixed model in \cite{mhk} is a variant of the HK model. For the mixed model, each agent can choose its degree of stubbornness and mix its opinion with the average opinion of it neighbors at each update. Let $[n]=\{1,2,\ldots,n\}$. The mixed model is as follows:
\begin{equation}\label{mHK}
x_i(t+1)=\alpha_i(t)x_i(t)+\frac{1-\alpha_i(t)}{|N_i(t)|}\sum_{j\in N_i(t)}x_j(t)    
\end{equation}
where
$$\begin{array}{rcl}
    \displaystyle x_i(t)&\n=\n& \hbox{opinion of agent $i$ at time $t$,}\vspace{4pt}\\
    \displaystyle N_i(t)&\n =\n& \{j\in [n]: \norm{x_i(t)-x_j(t)}\leq\epsilon\}\hbox{ is the neighborhood of agent $i$ at time $t$},\vspace{4pt}  \\
   \displaystyle \alpha_i(t)&\n\in \n&[0,1] \hbox{ is the degree of stubbornness of agent $i$} 
\end{array}$$
for all $i\in[n]$ and $t\in\mathbf{N}$. Written in matrix form,
$$x(t+1)=\big(\diag(\alpha(t))+(I-\diag(\alpha(t)))A(t)\big)x(t)$$ where
$$\begin{array}{rcl}
     \displaystyle x(t) &\n=\n& (x_1(t),\ldots,x_n(t))' = \hbox{transpose of }(x_1(t),\ldots,x_n(t)),\vspace{4pt}  \\
     \displaystyle \alpha(t) &\n =\n& (\alpha_1(t),\ldots,\alpha_n(t))' = \hbox{transpose of }(\alpha_1(t),\ldots,\alpha_n(t)),\vspace{4pt}\\
     \displaystyle A_{ij}(t) &\n =\n& \mathbbm{1}\{j\in N_i(t)\}/|N_i(t)|.
\end{array}$$
In particular, \eqref{mHK} reduces to
\begin{itemize}
    \item the synchronous HK model if $\alpha(t)=\Vec{0}$ for all $t\geq 0$, and\vspace{4pt}
    \item the asynchronous HK model if $\alpha(t)=(\mathbbm{1}\{k\neq i(t)\})_{k=1}^n$ for all $t\geq0$ and for some uniformly randomly selected $i(t)\in[n]$.
\end{itemize}
Say agent $i$ at time $t$ is
\begin{itemize}
    \item absolutely stubborn if $\alpha_i(t)=1$,\vspace{2pt}
    \item not absolutely stubborn, or open-minded if $\alpha_i(t)<1$, and\vspace{2pt}
    \item absolutely open-minded if $\alpha_i(t)=0$.
\end{itemize} 
In \cite{mhk}, we study the mixed model deterministically and point out some properties of the synchronous HK model, such as finite-time convergence, do not hold for the mixed model. Here, we assume that $$U_t=\{i\in [n]:\alpha_i(t)<1\}\hbox{ for all }t\geq0$$ are independent and identically distributed random variables with a support $S\subset\mathscr{P}([n])$ containing a partition of $[n]$, say $\{K_i\}_{i=1}^s$, where $\mathscr{P}([n])$ is the power set of $[n]$. Let $(\Omega,\mathscr{F},P)$ be a probability space for $\mathscr{F}\subset\mathscr{P}(\Omega)$ a $\sigma$-algebra and $P$ a probability measure. Observe that the asynchronous HK model is a particular case where the support $S=\{\{i\}\}_{i=1}^n$ is also a partition of $[n]$ and $U_t$ are uniform random variables on $S$, denoted by $U_t=\unif{(S)}$, for all $t\geq0$. 

\begin{definition}
A \emph{profile} at time $t$ is an undirected graph $\mathscr{G}(t)$ with the vertex set and edge set
$$\mathscr{V}(t)=[n]\quad \hbox{and}\quad \mathscr{E}(t)=\{ij:i\neq j\hbox{ and }\norm{x_i(t)-x_j(t)}\leq\epsilon\}.$$
\end{definition}

\begin{definition}
A profile $\mathscr{G}(t)$ is $\delta$-trivial if any two of its vertices are at a distance of at most $\delta$ apart.
\end{definition}

Denote almost surely by a.s., which may be omitted for simplicity. The theorems and corollary in \cite{mhk} can be interpreted nondeterministically as follows.
\begin{theorem}[\cite{mhk}]\label{t1}
     Let $$\beta_t := \max_{i, j \in [n], \alpha_i (t) \geq \alpha_j (t)} \bigg(\alpha_i (t) - \frac{\alpha_i (t) - \alpha_j (t)}{n} \bigg).$$
 Assume that~$\limsup_{t \to \infty} \beta_t < 1$ a.s. and that~$\mathscr{G} (t)$ is~$\epsilon$-trivial a.s.. Then,
 $$ \lim_{t \to \infty} \max_{i, j \in [n]} \|x_i (t) - x_j (t) \| = 0\hbox{ a.s.}. $$
\end{theorem}
For an $\epsilon$-trivial profile, agents need not play open-minded all the time. As long as there are infinitely many $\beta_t$ having an upper bound less than 1, eventually can a consensus be achieved. Given an $\epsilon$-trivial profile of the synchronous HK model, it is clear that a consensus can be achieved at the next time step. Observe that $\limsup_{t\to\infty}\beta_t<1$ automatically holds for the asynchronous HK model; therefore a consensus can be achieved eventually given an $\epsilon$-trivial profile.

\begin{theorem}[\cite{mhk}]\label{t2}
 Define~$d_t^i = \max_{j \in N_i (t)} \|x_i (t) - x_j (t) \|$. If
 $$ \sum_{t = 0}^{\infty} \ (1 - \alpha_i (t)) \bigg(1 - \frac{1}{|N_i (t)|} \bigg) d_t^i < \infty\hbox{ a.s.}, \ \hbox{then} \ x_i (t) \to x_i \in \mathbf{R^d} \ \hbox{as} \ t \to \infty\hbox{ a.s..} $$ 
\end{theorem}
It is difficult to track the dynamics trajectories. Nevertheless, $\alpha_i$ is controllable and $\big(1 - \frac{1}{|N_i (t)|} \big) d_t^i$ is bounded. Therefore, the assumption of Theorem \ref{t2} holds as long as $\sum_{t\geq0}(1-\alpha_i(t))<\infty.$ Namely, an agent can achieve its opinion limit without considering the others.

\begin{theorem}[\cite{mhk}]\label{t3}
 Assume that~$\limsup_{t \to \infty} \max_{i \in [n]} \alpha_i (t) < 1$ a.s..
 Then, for any~$\delta > 0$, every component of a profile is~$\delta$-trivial in finite time a.s., i.e.,
 $$ \tau_{\delta} := \inf \{t \geq 0 : \hbox{every component of $\mathscr{G} (t)$ is $\delta$-trivial} \} < \infty\hbox{ a.s.}. $$
\end{theorem}
All components of a profile is $\delta$-trivial in finite time as long as all agents are open-minded infinitely many times and their degree of stubbornness has an upper bound less than 1.

\begin{corollary}[\cite{mhk}]
\label{co1}
 Assume that~$\sup_{t \in \mathbf{N}} \max_{i \in [n]} \alpha_i (t) < 1$ a.s..
 Then, $\tau_{\delta}$ is bounded from above a.s..
 Also, letting~$\hat{\tau}_m = \tau_{\epsilon / m}$ for~$m \geq 4$, there is no interactions between any two components of~$\mathscr{G} (t)$ at the next time step
 for some $M \geq 4$ and for all~$t \geq \hat{\tau}_M$ a.s., i.e.,
 $$ \mathscr{G} (t) = \mathscr{G} (\hat{\tau}_M) \quad \hbox{for some} \quad M \geq 4 \quad \hbox{and for all} \quad t \geq \hat{\tau}_M\hbox{ a.s.}. $$
 In particular, $x$ in \eqref{mHK} is asymptotically stable a.s..
\end{corollary}
All components of a profile are $\delta$-trivial by some finite time given that all agents are open-minded and their degree of stubbornness has an upper bound less than 1 all the time. 

\section{main results}
For Theorems \ref{t1}, \ref{t3} and Corollary \ref{co1}, the assumptions involve the whole agents. However, the following results, such as asymptotic stability, hold without considering the absolutely stubborn who can vary over time. In particular, the following results can interpret the asynchronous HK model.

\begin{theorem}\label{T3}
Assume that $0\leq\limsup_{t\to\infty}\sup\{\alpha_i(t):i\in[n]\hbox{ and } \alpha_i(t)<1\}<\gamma<1$ for some $\gamma$ constant a.s.. Then, for any $\delta>0,$ all components of a profile are $\delta$-trivial in finite time a.s., i.e.,
 $$ \tau_{\delta} := \inf \{t \geq 0 : \hbox{all components of $\mathscr{G} (t)$ is $\delta$-trivial} \} < \infty\hbox{ a.s.}. $$

\end{theorem}
All components of a profile are $\delta$-trivial in finite time given that there are open-minded agents infinitely many times and their degree of stubbornness has a constant upper bound less than 1.

\begin{corollary}
\label{Co1}
 Assume that~$0\leq \sup \{\alpha_i (t):i\in[n]\hbox{ and }\alpha_i(t)<1\}\leq\gamma< 1$ for some $\gamma$ constant a.s..
 Then, $E(\tau_{\delta})$ is bounded from above.  Also, letting~$\hat{\tau}_m = \tau_{\epsilon / m}$ for~$m \geq 4$, there is no interactions between any two components of~$\mathscr{G} (t)$ at the next time step
 for some $M \geq 4$ and for all~$t \geq \hat{\tau}_M$ a.s., i.e.,
 $$ \mathscr{G} (t) = \mathscr{G} (\hat{\tau}_M) \quad \hbox{for some} \quad M \geq 4 \quad \hbox{and for all} \quad t \geq \hat{\tau}_M\hbox{ a.s.}. $$
 In particular, $x$ in \eqref{mHK} is asymptotically stable a.s..
 
\end{corollary}
The expected number of time steps until all components of a profile is $\delta$-trivial is bounded from above. It turns out that for nondeterministic mixed model, still can the results, such as asymptotic stability, hold without considering the absolutely stubborn. 

\section{the mixed model}

\begin{definition}
 The \emph{termination time} of~$n$ agents, $T_n$, is the maximum number of iterations in~\eqref{mHK} by reaching a steady state over all initial profiles, i.e.,
 $$ T_n = \inf \{t\geq0 : x (t) = x (s) \ \hbox{for all} \ s \geq t \}. $$
\end{definition}

\begin{definition}
 A \emph{merging time} is a time~$t$ that two agents with different opinions at time~$t - 1$ have the same opinion at time~$t$, i.e.,
 $$ x_i (t) = x_j (t) \quad \hbox{and} \quad x_i (t - 1) \neq x_j (t - 1) \quad \hbox{for some} \quad i, j \in [n]. $$
\end{definition}

\begin{definition}
 The \emph{convex hull} generated by~$v_1, v_2, \ldots, v_n \in \mathbf{R^d}$ is the smallest convex set containing~$v_1, v_2, \ldots, v_n$, i.e.,
 $$ C (\{v_1, v_2 \ldots, v_n \}) = \{v : v = \sum_{i = 1}^n \lambda_i v_i \ \hbox{where} \ (\lambda_i)_{i = 1}^n \ \hbox{is stochastic} \}. $$
\end{definition}

\begin{definition}
 For~$\delta > 0$, $x(t)$ in~\eqref{mHK} is a \emph{$\mathbf{\delta}$-equilibrium} if there is a partition
 $$ \{G_1, G_2, \ldots, G_m \} \ \hbox{of the set} \ \{x_1 (t), x_2 (t), \ldots, x_n (t) \} $$
 such that the following two conditions hold:
 $$ \dist (C (G_i), C (G_j)) > \epsilon \ \hbox{for all} \ i \neq j \quad \hbox{and} \quad \diam (C (G_i)) \leq \delta \ \hbox{for all} \ i \in [m]. $$
\end{definition}
\cite{mhk} lists three properties distinct from the synchronous HK model.
\begin{enumerate}
    \item The termination time is not finite.\vspace{4pt}
    \item Agents merging at time~$t$ may depart at time~$t + 1$.
 In particular, $\mathscr{G}(t)$ $\epsilon$-trivial may not imply that~$x (t + 1)$ in~\eqref{mHK} is a steady state.\vspace{4pt}
 \item A $\delta$-equilibrium may not exist for all $0<\delta \leq \epsilon$.
\end{enumerate}
The following are properties differing from the asynchronous HK model.
\begin{property}
Merging can exist.
\end{property}

\begin{example}
Consider $x_1(0)=0,\  x_2(0)=\epsilon,\ \alpha_1(0)=0\hbox{ and }\alpha_2(0)=0$. Then, $x_1$ and $x_2$ merge at time $t=1$.
\end{example}

\begin{property}
A $\delta$-equilibrium may not exist for all $0<\delta\leq\epsilon.$
\end{property}

\begin{example}
See Example 3 in \cite{mhk}.
\end{example}

\begin{lemma}[\cite{mhk}]\label{NL8}
 Let~$Z(t)=\sum_{i,j\in[n]}\|x_i(t)-x_j(t)\|^2\wedge\epsilon^2.$ Then,~$Z$ is nonincreasing with respect to~$t$. In particular,
$$\begin{array}{rcl}
    \displaystyle Z(t)-Z(t+1)&\n \geq\n &\displaystyle 4 \sum_{i=1}^n\bigg(1+|N_i(t)|\frac{\alpha_i(t)}{1-\alpha_i(t)}\mathbbm{1}\{\alpha_i(t)<1\}\bigg)\|x_i(t)-x_i(t+1)\|^2.
\end{array}$$
\end{lemma}

\begin{definition}{\rm
 A symmetric matrix~$M$ is called a \emph{generalized Laplacian} of a graph $G = (V, E)$ if for~$x, y \in V$, the following two conditions hold:
 $$ M_{xy} = 0 \ \hbox{for} \  x \neq y \ \hbox{and} \ xy \notin E \quad \hbox{and} \quad M_{xy} < 0 \ \hbox{for} \ x \neq y \ \hbox{and} \ xy \in E. $$
 Let~$d_G (x)$ = degree of~$x$ in $G$, let~$V (G)$ = vertex set of~$G$, and let~$E (G)$ = edge set of~$G$.
 Then, the \emph{Laplacian} of~$G$ is defined as~$\mathscr{L} = D_G - A_G$ where
 $$ D_G = \diag ((d_G (x))_{x \in V (G)}) \quad \hbox{and} \quad A_G = \ \hbox{the adjacency matrix}. $$
 In particular, $(A_G)_{xy} = \mathbbm{1} \{xy \in E (G)\}$ when the graph~$G$ is simple.}
\end{definition}

 Observe that there is no restrictions on the diagonal entries of the matrix~$M$.
 Moreover, the Laplacian of~$G$ is a generalized Laplacian.

\begin{proof}[\bf Proof of Theorem \ref{T3}]
By the assumption, there is $(t_k)_{k\geq1}\subset\mathbf{N}$ strictly increasing such that $\max\{\alpha_i(t_k):i\in[n]\hbox{ and }\alpha_i(t_k)<1\}<\gamma_1$ for some $\gamma<\gamma_1<1$ constant and for all $k\geq1$.
For all $m\geq 1$,
\begin{equation}\label{tel}
    n^2\epsilon^2>Z(0)\geq Z(0)-Z(m)=\sum_{t=0}^{m-1}[Z(t)-Z(t+1)].
\end{equation}
Letting $m\to\infty$ and applying Lemma \ref{NL8},
\begin{equation*}
    n^2\epsilon^2\geq\sum_{t\geq 0}[Z(t)-Z(t+1)]\geq 4\sum_{t\geq 0}\sum_{i\in[n],\alpha_i(t)<1}\norm{x_i(t)-x_i(t+1)}^2.
\end{equation*}
Set $\tau=\tau_{\delta}$. Assume by contradiction that $\tau=\infty$ on some $E\in \mathscr{F}$ with $P(E)>0$. Taking expectation on both sides on $E$, denoted by $E_E$,
\begin{align*}
    n^2\epsilon^2\ &\geq\ 4\sum_{t\geq0}\sum_{j=1}^s E_E\bigg(\sum_{i\in[n],\alpha_i(t)<1}\norm{x_i(t)-x_i(t+1)}^2\bigg|U_t=K_j\bigg)P(U_t=K_j)\\
    &=\ 4\sum_{t\geq 0}\sum_{j=1}^s P(U_t=K_j)E_E\bigg(\sum_{i\in K_j}\norm{x_i(t)-x_i(t+1)}^2\bigg)\\
    &\geq\ 4\min_{j\in[s]}P(U_0=K_j)\sum_{t\geq 0}E_E\bigg(\sum_{j=1}^s \sum_{i\in K_j}\norm{x_i(t)-x_i(t+1)}^2\bigg)\\
    &=\ 4\min_{j\in[s]}P(U_0=K_j)\sum_{t\geq 0}E_E\bigg(\sum_{i=1}^n \norm{x_i(t)-x_i(t+1)}^2\bigg)\\
    &>\ 4\min_{j\in[s]}P(U_0=K_j)\sum_{t\geq 0}E_E\bigg(\frac{2\delta^2(1-\max_{i\in [n]}\alpha_i(t))^2}{n^8}\bigg) \tag{+}\label{plu}\\
    &\geq\ 4\min_{j\in[s]}P(U_0=K_j)E_E\bigg(\sum_{k\geq 1}\frac{2\delta^2(1-\max_{i\in [n]}\alpha_i(t_k))^2}{n^8}\bigg)\\
    &\geq\ 4\min_{j\in[s]}P(U_0=K_j)\sum_{k\geq 1}E_E\bigg(\frac{2\delta^2(1-\gamma_1)^2}{n^8}\bigg)=\infty,\hbox{ a contradiction,}
\end{align*}
where \eqref{plu} refers to the proof of Theorem 3 in \cite{mhk}. Therefore, $\tau$ is finite. 

\end{proof}

\begin{proof}[\bf Proof of Corollary \ref{Co1}] 
Since the assumption of Corollary \ref{Co1} meets that of Theorem \ref{T3}, $\tau$ is finite a.s.. Setting $m=\tau$ in \eqref{tel} and applying Lemma \ref{NL8}, $$n^2\epsilon^2>\sum_{t=0}^{\tau-1}[Z(t)-Z(t+1)]\geq4\sum_{t=0}^{\tau-1}\sum_{i\in[n],\alpha_i(t)<1}\norm{x_i(t)-x_i(t+1)}^2.$$
Similarly as \eqref{plu},
\begin{align*}
    n^2\epsilon^2>&\ 4\min_{j\in[s]}P(U_0=K_j)E\bigg(\sum_{t=0}^{\tau-1}\frac{2\delta^2(1-\max_{i\in [n]}\alpha_i(t))^2}{n^8}\bigg)\\
    \geq &\ 4\min_{j\in[s]}P(U_0=K_j)E\bigg(\frac{2\tau\delta^2(1-\gamma)^2}{n^8}\bigg)
\end{align*}
therefore $$E(\tau)<\frac{n^{10}}{8(1-\gamma)^2\min_{j\in[s]}P(U_0=K_j)}(\frac{\epsilon}{\delta})^2.$$
Hence, $E(\tau)$ is bounded from above. 

To show the asymptotic stability of~$x$, first observe that~$\hat{\tau}_m$ is finite and nondecreasing with respect to~$m$. For all~$0<\delta\leq\epsilon/4$ and~$t\geq 0$, the following conditions are equivalent given that every component of~$\mathscr{G}(t)$ is~$\delta$-trivial:
 \begin{enumerate}
     \item\label{q1} Some component of~$\mathscr{G}(t+1)$ is~$\delta$-nontrivial.\vspace*{2pt}
     \item\label{q2} Some components of~$\mathscr{G}(t)$ interact at time~$t+1$.\vspace*{2pt}
     \item\label{q3} Some component of~$\mathscr{G}(t+1)$ is~$\epsilon/2$-nontrivial.
 \end{enumerate}
 It is clear that~\ref{q1} $\Rightarrow$ \ref{q2} and~\ref{q3} $\Rightarrow$ \ref{q1}; therefore we show~\ref{q2} $\Rightarrow$ \ref{q3}.
 \begin{proof}[Proof of \ref{q2} $\Rightarrow$ \ref{q3}]\noindent
 Let the convex hull of a component~$G$ be
 $$ Cv(G) = C(\{x_j : j \in V(G) \}). $$
 The fact that some components of~$\mathscr{G} (t)$ interact at time~$t+1$ implies that there exist
 $$ i, j \in [n] \quad \hbox{with} \quad ij \in \mathscr{E} (t + 1), \quad \quad i\in V(G_{\Tilde{i}})\quad \hbox{and}\ j\in V(G_{\Tilde{j}}) $$
 for some distinct components~$G_{\Tilde{i}}$ and~$G_{\Tilde{j}}$ of~$\mathscr{G} (t)$.
 Therefore,
 $$ x_i (t + 1) \in Cv (G_{\Tilde{i}}) \quad \hbox{and} \quad x_j (t + 1) \in Cv (G_{\Tilde{j}}). $$
 Hence,
\begin{align*}
 \epsilon \ < \ & \|x_i (t) - x_j (t) \| \\
          \ \leq \ & \|x_i (t) - x_i (t + 1) \| + \|x_i (t + 1) - x_j (t + 1) \| + \|x_j (t + 1) - x_j (t) \| \\
          \ \leq \ & \delta + \|x_i (t + 1) - x_j (t + 1) \| + \delta = \|x_i (t + 1) - x_j (t + 1) \| + 2 \delta. \end{align*}
This implies that
$$ \|x_i (t + 1) - x_j (t + 1) \| > \epsilon - 2 \delta \geq \epsilon - 2 \cdot \frac{\epsilon}{4} = \frac{\epsilon}{2} \quad \hbox{for all} \quad 0 < \delta \leq \frac{\epsilon}{4} $$
so the component of~$\mathscr{G} (t + 1)$ containing~$ij$ is~$\epsilon/2$-nontrivial.
 \end{proof}\noindent
 Let
 $$ A_m = \{t \in [\hat{\tau}_m, \hat{\tau}_{m + 1}) : \hbox{some component of $\mathscr{G}(t)$ is $\epsilon/m$-nontrivial} \} $$
 and~$t_m = \inf A_m$.

 $$\hbox{Claim:}\quad \hbox{the set}\ \mathscr{A} := \{t_k : A_k \neq \varnothing \} \  \hbox{is finite a.s.}.$$
For~$t_m \in \mathscr{A}$, since some component of~$\mathscr{G} (t_m)$ is~$\epsilon/m$-nontrivial and all components of~$\mathscr{G} (t_m - 1)$
 are~$\epsilon/m$-trivial, by~\ref{q1}$\Rightarrow$ \ref{q3}, some component of~$\mathscr{G} (t_m)$ is~$\epsilon/2$-nontrivial.
 Similarly as \eqref{plu}, we get
\begin{align*}
  n^2 \epsilon^2 \ \geq \ & 4\min_{j\in[s]}P(U_0=K_j)E\bigg( \sum_{t \in\mathscr{A}} \frac{2(\epsilon/2)^2(1-\max_{i\in[n]}\alpha_i(t))^2}{n^8}\bigg)\\
  \geq\ & 4\min_{j\in[s]}P(U_0=K_j)E\bigg( \sum_{t \in\mathscr{A}} \frac{2(\epsilon/2)^2(1-\gamma)^2}{n^8}\bigg)\\
  =\ & E(|\mathscr{A}|)\frac{2\epsilon^2(1-\gamma)^2}{n^8}\min_{j\in[s]}P(U_0=K_j), 
  \end{align*}
 from which it follows that
 $$E( |\mathscr{A}| )\leq \frac{n^{10}}{2(1 - \gamma)^2\min_{j\in[s]}P(U_0=K_j)}. $$
 Hence, the set~$\mathscr{A}$ is finite a.s..
 By the fact that~$\mathscr{A}$ is finite and that~\ref{q2}$\Rightarrow$ \ref{q1}, there is no interactions between any two components of~$\mathscr{G} (s)$ at the next time step for some~$M \geq 4$ and for all~$s \geq \hat{\tau}_M$. Hence, all components of~$\mathscr{G} (\hat{\tau}_M)$ are independent systems. Therefore, $x$ in~\eqref{mHK} is asymptotically stable.

\end{proof}
Note that Corollary \ref{Co1} particularly verifies that the asymptotic stability holds for the asynchronous HK model. Theorem 5 in \cite{ts} states that the expected number of steps until all agents in the asynchronous HK model reach a $\delta$-equilibrium is bounded from above. Therefore, asymptotic stability holds for the asynchronous HK model. In fact, a $\delta$-equilibrium is not necessary for asymptotic stability.  Since the assumption of Theorem \ref{t3} automatically holds for the asynchronous HK model, there is no interactions between any two components of $\mathscr{G}(t)$ at the next time step for some $M\geq4$ and for all $t\geq\hat{\tau}_M$ a.s.; therefore asymptotic stability holds. Also, the assumption of Corollary \ref{Co1} automatically holds for the asynchronous HK model; therefore $E(\hat{\tau}_{M})$ is bounded from above.

\section{conclusion}
The mixed model covers both the synchronous HK model and the asynchronous HK model. Each agent can decide its degree of stubbornness and mix its opinion with the average opinion of its neighbors at all times. Compared to the HK dynamics in which an agent is either absolutely stubborn or absolutely open-minded, agents for the mixed model are more flexible and closer to reality, making it harder to attain asymptotic stability. We elaborate in \cite{mhk} that the synchronous HK model is in fact a particular case of the deterministic mixed model. In this theme, we illustrate that the asynchronous model is in fact a particular case of the nondeterministic mixed model. Furthermore, under some circumstances the asymptotic stability can hold without considering the absolutely stubborn for the nondeterministic mixed model.


\begin{thebibliography}{99}
\bibitem{mhk}
\newblock H. Li,
\newblock Mixed Hegselmann-Krause Dynamics,
\newblock \emph{Discrete and Continuous Dynamical Systems-B}, (2021), \url{https://arxiv.org/abs/2010.03050}

\bibitem{ts} %(MR3365075) [10.1109/TAC.2015.2394954]
\newblock S. R. Etesami and T. Basar,
\newblock Game-theoretic analysis of the Hegselmann-Krause model for Opinion dynamics in finite dimensions,
\newblock \emph{IEEE Transactions on Automatic Control}, \textbf{60} (2015), 1886--1897.

\bibitem{6}
\newblock L. W. Beineke, P. J. Cameron and R. J. Wilson,
\newblock \emph{Topics in Algebraic Graph Theory},
\newblock Cambridge University Press, Cambridge, UK, 2004.

\bibitem{5} %(MR2340484) [10.1007/978-3-540-73510-6]
\newblock T. Biyikoglu, J. Leydold and P. F. Stadler,
\newblock \emph{Laplacian Eigenvectors of Graphs: Perron-Frobenius and Faber-Krahn Type Theorems,}
\newblock Springer-Verlag, Berlin Heidelberg, 2007.




\bibitem{7} %(MR2978290)
\newblock R. A. Horn and C. R. Johnson,
\newblock \emph{Matrix Analysis,}
\newblock Cambridge University Press, Cambridge, 2013.



\end{thebibliography}
\end{document}